\UseRawInputEncoding
\documentclass[reqno,11pt]{amsart}
\usepackage{amsthm,amssymb,latexsym,mathrsfs}
\usepackage{amsmath}
\usepackage{indentfirst}
\usepackage{color}

\usepackage[top=2.5cm,bottom=2.5cm,left=2.5cm,right=2.5cm]{geometry}
\allowdisplaybreaks[1]
\numberwithin{equation}{section}
\newtheorem{theorem}{Theorem}[section]

\newtheorem{lemma}[theorem]{Lemma}
\newtheorem{remark}[theorem]{Remark}

\newtheorem{corollary}[theorem]{Corollary}
\newtheorem{assumption}[theorem]{Assumption}
\allowdisplaybreaks

\begin{document}

\title[Mean-field limit]{Convergence towards the population cross-diffusion system from stochastic many-particle system}

\author{Yue Li}
\address{Department of Mathematics, Nanjing University,
 Nanjing 210093, P.R. China}
\email{liyue2011008@163.com}

\author{Li Chen}
\address{Lehrstuhl f\"{u}r Mathematik IV, Universit\"{a}t Mannheim,
Mannheim 68131, Germany}
\email{chen@math.uni-mannheim.de}

\author{Zhipeng Zhang}
\address{Department of Mathematics, Nanjing University,
 Nanjing 210093, P.R. China}
\email{zhangzhipeng@nju.edu.cn}

%




\begin{abstract}
In this paper, we derive rigorously a non-local cross-diffusion system from an interacting stochastic many-particle system in the whole space. The convergence is proved in the sense of probability by introducing an intermediate particle system with a mollified interaction potential, where the mollification is of algebraic scaling. The main idea of the proof is to study the time evolution of a stopped process and obtain a Gr\"onwall type estimate by using Taylor's expansion around the limiting stochastic process.
\end{abstract}

\maketitle

\section{ Introduction}

In this paper, we give a rigorous justification for the mean-field limit from an interacting particle system to the population cross-diffusion system
as the number of particles goes to infinity.
More precisely, we present the derivation of $n$-species cross-diffusion system as follows
\begin{equation}\label{VNS}
\left\{
\begin{aligned}
  &\partial_t u_i={\rm{div}}(u_i\nabla U_i)+\sigma_i \Delta u_i+{\rm{div}}\Big[u_i\sum^n_{j=1}\nabla f(B_{ij}\ast u_j)\Big],\;
  B_{ij}(|x|)=\frac{C(d,\vartheta_{ij})}{|x|^{\vartheta_{ij}}},\;\vartheta_{ij}\in(0,d-2],\\
  &u_i(0)=u^0_i(x),\qquad i=1,\ldots,n,
\end{aligned}
\right.
\end{equation}
where $\sigma_i>0$ are the constant diffusion coefficients,
$\mathbf{u}=(u_1,\ldots,u_n)$ stands for the vector of population densities,
$U_i(x)=-\frac{1}{2}|x|^2$ represent environment potentials
and $C(d,\vartheta_{ij})$ are constants depend on $d$ and $\vartheta_{ij}$.
The transitions rates depend on the densities by a nonlinear term $f$.

The aim of this paper is to rigorously derive the system \eqref{VNS} from the following stochastic many-particle system. This system describes the movements of $n$ species of particles, with the particle numbers $N_i\in \mathbb{N}$ $(i=1,\ldots,n)$, according to the given law. Without loss of generality, we let $N=N_i$ ($i=1,\ldots,n$). Let $(\Omega, \mathcal{F}, (\mathcal{F}_{t\geq 0}) , \mathbb{P})$ be a complete filtered probability space.  We consider  $d$-dimensional $\mathcal{F}_t$-Brownian motions $(W^k_i(t))_{t\geq 0} $ ($k=1,\ldots,N$, $i=1,\ldots,n$) which are assumed to be independent of each other.  We assume that $(\xi^k_i)$ ($k=1,\ldots,N$, $i=1,\ldots,n$) are i.i.d. random variables, independent of $(W^k_i(t))_{t\geq 0}$, and have common probability density function $u_i^0$. We use the notation $X^{N,k}_{\eta,i}(t)$ to represent the $k$-th particle of $i$-th species and the dynamics of
$X^{N,k}_{\eta,i}(t)$ are governed by
\begin{equation}\label{1}
\left\{
\begin{aligned}
  &dX^{N,k}_{\eta,i}=\Big[-\nabla U_i(X^{N,k}_{\eta,i})-\sum^n_{j=1}\nabla f_\gamma\Big(\frac{1}{N}\sum^N_{l=1}B^\eta_{ij}(X^{N,k}_{\eta,i}-X^{N,l}_{\eta,j})\Big)\Big]dt
  +\sqrt{2\sigma_i}dW^k_i(t),\\
  &X^{N,k}_{\eta,i}(0)=\xi_i^k,\qquad i=1,\ldots,n, \qquad k=1,\ldots,N,
\end{aligned}
\right.
\end{equation}
where $f_\gamma$ is an approximation of $f$ which can be constructed, for example in Remark 1.2, and
\begin{align*}
	B^\eta_{ij} := \begin{cases}
		V^\eta\ast B_{ij}, \ &0<\vartheta_{i,j}<d-2,  \\
		V^\eta\ast\bar B_{ij}, \ &\vartheta_{i,j}=d-2,
	\end{cases}
\quad
\bar B_{ij} (|x|) := \begin{cases}
B_{ij}(|x|), \ &|x|\geq \eta,  \\
B_{ij}(\eta), \ &|x|< \eta.
\end{cases}
\end{align*}
Here $V^\eta(x):=\frac{1}{\eta^d}V(\frac{x}{\eta})$ with $\eta>0$ is a mollification kernel which means $V\geq 0$ is a given radially symmetric smooth
function such that $\int_{\mathbb{R}^d}V(x)\,dx=1$.

The problem considered in this paper dedicates to the understanding of diffusion (and cross-diffusion) effects on the microscopic level. It belongs to the research of mean-field limit for interacting particle system.
There have been extensive studies of the mean-field limit problems in the last decades. Many important contributions have been made for problems with singular interacting potentials such as the Coulomb potential in Keller-Segel systems. An extensive review of this research field is out of the scope of this paper, we refer to \cite{G,JW,S2020,CGHL} for more detailed summary on the results and methods.

The convergence of moderate interacting system was introduced and proved by Oelschl\"{a}ger in \cite{O1984,O1989,O1990} in order to derive reaction-diffusion equations and the porous medium equation.
The authors \cite{JM} considered further the fluctuation of this problem.
This idea has been used to derive chemotaxis equation from an interacting stochastic many-particle system in \cite{St}.
The derivation of cross-diffusion type systems has only been studied in the last few years.
It is proved in \cite{Seo} that the hydrodynamic limit of the empirical densities of two types is the solution to the Maxwell-Stefan equation. The authors in \cite{FM} derived the non-local Lotka-Volterra system with cross-diffusion from particle system.
The Shigesada-Kawasaki-Teramoto system was obtained from a microscopic many-particle Markov process in \cite{D}.
Rigorous derivation of the degenerate parabolic-elliptic Keller-Segel system from a moderate interacting stochastic particle system was given in \cite{CGHL}.
There are very few results for more than two species.
In \cite{CDJ1}, the authors established the global existence of weak solutions to cross-diffusion systems for an arbitrary number of competing population species.
The mean-field limit of a moderate interacting stochastic many-particle system for multiple population species is obtained in \cite{CDJ2} with the logarithmic scaling. Furthermore, with the same scaling the authors in \cite{CDHJ} derived population cross-diffusion systems of Shigesada-Kawasaki-Teramoto type from stochastic moderately interacting many-particle systems for multiple population species.
This paper is aimed to derive \eqref{VNS} from \eqref{1} with algebraic scaling $\eta=N^{-\beta}$ for some $\beta$.

To precisely state the main results of this paper, we first give the general assumptions on $f$.
\begin{assumption} \label{1.1}
	Let $f\in C^3(\mathbb{R}_*;[0,\infty))$, where $\mathbb{R}_*=\mathbb{R}_+\cup \{0\}$, and for some $m>0$ it holds
	\begin{align*}
	\|f\|_{C^3([-M,M])}\leq M^m {\quad \rm{for\;\; any}}\; M \gg 1.
	\end{align*}
\end{assumption}
\begin{remark}\label{fgamma}
	A possible approximation $f_\gamma\in C^3(\mathbb{R}_*;[0,\infty))$ of $f$ can be given by
\begin{align*}
	f_\gamma(r) := \begin{cases}
		f(r), \ & 0\leq r\leq\frac{1}{2\gamma},  \\
		f(\frac{1}{\gamma}), \ & \frac{1}{\gamma}\leq r.
	\end{cases}
\end{align*}
It is obviously that the following estimate holds
	\begin{align*}
	\|f_\gamma\|_{C^3(\mathbb{R})}\leq \|f\|_{C^3([-\frac{1}{\gamma},\frac{1}{\gamma}])}\leq \frac{1}{\gamma^m}.
	\end{align*}
\end{remark}

In order to prove the limit from \eqref{1} to \eqref{VNS}, we introduce an intermediate particle problem. This problem is formally viewed as a mean field limit  $N\rightarrow \infty$ in the system \eqref{1} for fixed $\eta,\gamma>0$, namely
\begin{equation}\label{2}
\left\{
\begin{aligned}
&d\bar X^{k}_{\eta,i}=\Big[-\nabla U_i(\bar X^{k}_{\eta,i})-\sum^n_{j=1}\nabla f_\gamma(B^\eta_{ij}\ast u_{\eta,j}(\bar X^{k}_{\eta,i}))\Big]dt
+\sqrt{2\sigma_i}dW^k_i(t),\\
&\bar X^{k}_{\eta,i}(0)=\xi_i^k,\qquad i=1,\ldots,n, \qquad k=1,\ldots,N.
\end{aligned}
\right.
\end{equation}
Here $u_{\eta,j}$ is the probability density function of $\bar X^{k}_{\eta,j}$ and satisfies the following cross-diffusion system:
\begin{equation}\label{3}
\left\{
\begin{aligned}
&\partial_t u_{\eta,i}={\rm{div}}(u_{\eta,i}\nabla U_i)+\sigma_i \Delta u_{\eta,i}
+{\rm{div}}\Big[u_{\eta,i}\sum^n_{j=1}\nabla f_\gamma(B^\eta_{ij}\ast u_{\eta,j})\Big],\\
&u_{\eta,i}(0)=u^0_i(x),\qquad i=1,\ldots,n.
\end{aligned}
\right.
\end{equation}

In this paper, we focus on the derivation of \eqref{VNS} from interacting particle system. To achieve this, we need the following assumptions of PDE solutions.
\begin{assumption}\label{as1.3}
Assume that $u_\eta,u\in L^\infty(0,T;L^1\cap H^s(\mathbb{R}^d))$ $(s>\frac{d}{2}+1)$ are solutions of systems \eqref{VNS} and \eqref{3}
respectively, furthermore it holds that
\begin{align}\label{est3}
\|u-u_\eta\|_{L^\infty(0,T;H^s(\mathbb{R}^d))}\leq C\eta,
\end{align}
where $C$ is a positive constant which is independent of $\eta$.
\end{assumption}
Actually, \eqref{VNS} and \eqref{3} are parabolic systems. The above assumption can be obtained at least for small initial data.
Similar to \cite{CDHJ}, the assumption \ref{as1.3} implies directly the wellposedness of SDE system \eqref{1} and the McKean-Vlasov problem of \eqref{VNS}, i.e.
\begin{equation}\label{4}
\left\{
\begin{aligned}
&d\hat X^{k}_{i}=\Big[-\nabla U_i(\hat X^{k}_{i})-\sum^n_{j=1}\nabla f(B_{ij}\ast u_{j}(\hat X^{k}_{i}))\Big]dt
+\sqrt{2\sigma_i}dW^k_i(t),\\
&\hat X^{k}_{i}(0)=\xi_i^k,\qquad i=1,\ldots,n, \qquad k=1,\ldots,N,
\end{aligned}
\right.
\end{equation}
where $u_i$ solves the limiting cross-diffusion system \eqref{VNS} and is the probability density function of $\hat X^{k}_{i}$.
Namely, when additional $\displaystyle\int_{\mathbb{R}^d}|x|^2u_0(x)dx<\infty$, there exist unique square-integrable adapted stochastic processes with continuous
paths, which are strong solutions to systems \eqref{2} and \eqref{4}, respectively.

Therefore, \eqref{est3} provides directly the following estimate:
\begin{align*}
&\mathbb{E}\Big(\sum^n_{i=1}\sup_{0\leq s\leq T}\max_{1\leq k\leq N}|\bar X^{k}_{\eta,i}(s)-\hat X^{k}_{i}(s)|^2 \Big)\\
\leq &C\mathbb{E}\Big(\sum^n_{i=1}\max_{1\leq k\leq N}\int_0^T |\nabla U_i(\bar X^k_{\eta,i})(t)-\nabla U_i(\hat X^{k}_{i})(t)|^2dt \Big)\\
&+C\mathbb{E}\Big(\sum^n_{i=1}\max_{1\leq k\leq N}\int_0^T\Big|\sum^n_{j=1}\big[\nabla f_\gamma(B^\eta_{ij}\ast u_{\eta,j}(\bar X^{k}_{\eta,i} ))
-\nabla f(B_{ij}\ast u_j(\hat X^{k}_{i}))\big] \Big|^2dt \Big)\\
\leq &C\int_0^T\mathbb{E}\Big(\sum^n_{i=1}\sup_{0\leq s\leq t}\max_{1\leq k\leq N}|\bar X^{k}_{\eta,i}(s)-\hat X^{k}_{i}(s)|^2 \Big)dt
+C\eta^{\frac{\vartheta_{ij}}{d}},
\end{align*}
where $C>0$ is a positive constant which is independent of $N$ and $\eta$. And combining Gr\"onwall inequality, we have
\begin{equation}
\label{est4}\mathbb{E}\Big(\displaystyle\sum^n_{i=1}\sup_{0\leq t\leq T}\max_{1\leq k\leq N}|\bar X^{k}_{\eta,i}(t)-\hat X^{k}_{i}(t)|^2\Big)
\leq C\eta^\frac{\vartheta_{ij}}{d}.
\end{equation}

The main result of this paper is the following
\begin{theorem}\label{main}
Let the assumptions \ref{1.1} and \ref{as1.3} hold, $0\leq u_0\in  L^1(\mathbb{R}^d)$, and $\displaystyle\int_{\mathbb{R}^d}|x|^2u_0(x)dx<\infty$. Assume that $T>0$, $\eta=N^{-\beta}$, $\gamma=N^{-\frac{\beta}{m}}$,
where $\beta\in \big(0,\frac{1}{2(5\sup_{1\leq i,j\leq n}\vartheta_{ij}+6)}\big)$, then for any arbitrary $\lambda>0$, it holds
\begin{align*}
\sup_{0\leq t\leq T}\mathbb{P}\Big(\sum^n_{i=1}\max_{1\leq k\leq N}|X^{N,k}_{\eta,i}(t)-\bar X^{k}_{\eta,i}(t)|>N^{-\alpha}\Big)\leq C(\lambda)N^{-\lambda},
\end{align*}
where $\alpha<\frac{1}{2}-\beta(2\sup_{1\leq i,j\leq n}\vartheta_{ij}+2)$ and $C(\lambda)$ is a positive constant independent of $N$.
\end{theorem}
Combined with the estimate in \eqref{est4}, we obtain the mean field limit result on the trajectory level and the propagation of chaos as a corollary
\begin{corollary} Under the same assumptions as in theorem \ref{main}, we have for any $\tilde\beta<\frac{\beta\vartheta_{ij}}{2d}$
\begin{align*}
\sup_{0\leq t\leq T}\mathbb{P}\Big(\sum^n_{i=1}\max_{1\leq k\leq N}|X^{N,k}_{\eta,i}(t)-\hat X^{k}_{i}(t)|>N^{-\tilde\beta}\Big)\leq CN^{-(\beta-\tilde\beta)}.
\end{align*}
Let $l\in \mathbb{N}$ and consider an $l$-tuple $(X^{N,1}_{\eta,i}(t), \ldots , X^{N,l}_{\eta,i}(t))$. We denote by  $P^{N,l}_{\eta,i}(t)$ its joint distribution.  Then it holds that
\begin{align*}
P^{N,l}_{\eta,i}(t) \ \text{converges weakly to} \ P_i^{\otimes l}(t)  \mbox{ as } N\rightarrow \infty,
\end{align*}
where $P_i(t)$ is a measure which is absolutely continuous with respect to the Lebesgue measure and has a probability density function $u_i(t,x)$.
\end{corollary}

The main result of this paper gives the same propagation of chaos result under the algebraic scaling $\eta\sim 1/N^\beta$. This result is obtained through the convergence in the sense of probability on the particle level. The benefit of algebraic scaling is that one can capture the singular interaction to some extend. To overcome the difficulty originated from the singular interaction, a suitable stopped process is established. Based on this, it is reduced to estimate the expectation of the stopped process with the help of Markov's inequality. A generalized version of Law of Large Numbers is the key point when we study the expectation. Another difficulty is caused by the nonlinear term. We have to find an approximate function $f_\gamma$ and give explicit scaling between $\gamma$ and $N$. Section 2 is dedicated to the proof of the main theorem.

\section{The proof of Theorem \ref{main}}\label{app}
We prove the convergence in probability on the particle level. Because of the singular interaction, one can not expect that under the algebraic scaling the convergence can be obtained in the expectation sense. The convergence in probability means that one allows that the particle trajectories are not always close, but the probability that they are not close is very low. Actually we can prove that the probability has a arbitrary convergence rate.

For any $\kappa\in\mathbb{N}$, we define a stopping time $\tau_\alpha$, a random variable $S^\kappa_\alpha$ and a set $B_\alpha$
\begin{align*}
&\tau_\alpha(\omega):=\inf\Big\{t\in(0,T):\sum^n_{i=1}\max_{1\leq k\leq N}|X^{N,k}_{\eta,i}(t)-\bar X^{k}_{\eta,i}(t)|\geq N^{-\alpha} \Big\}, \quad\omega\in\Omega,\\
&S^k_\alpha(t):= N^{\alpha \kappa}\sum^n_{i=1}\max_{1\leq k\leq N}|X^{N,k}_{\eta,i}(t\wedge\tau_\alpha)-\bar X^{k}_{\eta,i}(t\wedge\tau_\alpha)|^\kappa\leq 1,\\
&B_\alpha(t):=\{\omega\in\Omega:S^\kappa_\alpha(t)=1\}.
\end{align*}
By Markov's inequality, it holds
\begin{align*}
&\mathbb{P}\Big(\sum^n_{i=1}\max_{1\leq k\leq N}|X^{N,k}_{\eta,i}(t)-\bar X^k_{\eta,i}(t)|\geq N^{-\alpha} \Big)\\
\leq & \mathbb{P}\Big(\sum^n_{i=1}\max_{1\leq k\leq N}|X^{N,k}_{\eta,i}(t\wedge\tau_\alpha)-\bar X^k_{\eta,i}(t\wedge\tau_\alpha)|= N^{-\alpha} \Big)
\leq \mathbb{E}(S^\kappa_\alpha(t)).
\end{align*}
We notice that the introduction of parameter $\kappa$ is to increase the convergence rate. Actually, the above Markov's inequality works for arbitrary $\kappa$.
In order to complete the proof of Theorem \ref{main}, we just need to show that for any $\lambda>0$ and $T>0$,
it holds $\mathbb{E}(S^\kappa_\alpha(t))\leq CN^{-\lambda}$,
where the letter $C$ appeared in this section is a generic positive constant independent of $N$.
To this end, we need the following Law of Large Numbers, which can be found for example in \cite{CGHL}:
\begin{lemma}\label{large}
For $\varphi_{ij}\in L^\infty(\mathbb{R}^d)$, $i,j=1,\ldots,n$, we define for arbitrary $\theta\in (0,\frac{1}{2})$
\begin{align*}
A^{N,n}_{\theta,\varphi}(s):=\bigcup^n_{i,j=1}\bigcup^N_{k=1}
\Big\{\omega\in\Omega: \Big|\frac{1}{N}\sum^N_{l=1}\varphi_{ij}(\bar X^k_{\eta,i}(s)-\bar X^l_{\eta,j}(s))-\varphi_{ij}\ast u_{\eta,j}
(\bar X^k_{\eta,i}(s)) \Big|>\frac{1}{N^\theta} \Big\}.
\end{align*}
Then for any $m\in\mathbb{N}$, it holds that
\begin{align*}
P(A^{N,n}_{\theta,\varphi}(s))\leq C(n,m)N^{2m(\theta-\frac{1}{2})+1}
\big(\sup_{1\leq i,j\leq n}\|\varphi_{ij}\|^{2m}_{L^\infty(\mathbb{R}^d)}+\sup_{1\leq i,j\leq n}\|\varphi_{ij}\ast u_{\eta,j}\|^{2m}_{L^\infty((0,T)\times\mathbb{R}^d)} \big), \quad  \forall s\in [0,T].
\end{align*}
\end{lemma}
Next we study the time evolution of the cut-offed process $S^\kappa_\alpha(t)$. Notice that
\begin{align*}
&|X^{N,k}_{\eta,i}(t\wedge\tau_\alpha)-\bar X^{k}_{\eta,i}(t\wedge\tau_\alpha)|^\kappa\\
\leq & C\int_0^{t\wedge\tau_\alpha}|\nabla U_i(X^{N,k}_{\eta,i})-\nabla U_i(\bar X^{k}_{\eta,i}) |^\kappa\,ds\\
\qquad& +C\int_0^{t\wedge\tau_\alpha}\Big|\sum^n_{j=1}\nabla f_\gamma\Big(\frac{1}{N}\sum^N_{l=1}B^\eta_{ij}(X^{N,k}_{\eta,i}-X^{N,l}_{\eta,j})\Big)
  -\sum^n_{j=1}\nabla f_\gamma(B^\eta_{ij}\ast u_{\eta,j}(\bar X^{k}_{\eta,i}))   \Big|^\kappa\,ds\\
=:& J^1_{k,i}(t)+J^2_{k,i}(t).
\end{align*}
From the definition of $U_i(x)=-\frac{|x|^2}{2}$, we get
\begin{align}\label{5}
\mathbb{E}\Big(N^{\alpha \kappa}\sum^n_{i=1}\max_{1\leq k\leq N} J^1_{k,i}(t)\Big)
\leq
 &C\int_0^t \mathbb{E}(S^\kappa_\alpha(s))\,ds.
\end{align}
For the second term $J^2_{k,i}$,
\begin{align}\label{6}
&\mathbb{E}\Big(N^{\alpha \kappa}\sum^n_{i=1}\max_{1\leq k\leq N} J^2_{k,i}(t)\Big)\nonumber\\
\leq &C\mathbb{E}\Big(N^{\alpha\kappa}\sum^n_{i=1}\max_{1\leq k\leq N}\int_0^{t\wedge\tau_\alpha}
\Big|\sum^n_{j=1}f'_\gamma\Big(\frac{1}{N}\sum^N_{l=1}B^\eta_{ij}(X^{N,k}_{\eta,i}-X^{N,l}_{\eta,j})\Big)\nonumber\\
&\qquad\qquad\qquad\qquad\qquad\cdot\frac{1}{N}\sum^N_{l=1}\Big[\nabla B^\eta_{ij}(X^{N,k}_{\eta,i}-X^{N,l}_{\eta,j})
-\nabla B^\eta_{ij}(\bar X^{k}_{\eta,i}-\bar X^{l}_{\eta,j}) \Big] \Big|^\kappa\,ds \Big)\nonumber\\
&+C\mathbb{E}\Big(N^{\alpha\kappa}\sum^n_{i=1}\max_{1\leq k\leq N}\int_0^{t\wedge\tau_\alpha}\Big|\sum^n_{j=1}
\Big[f'_{\gamma}\Big(\frac{1}{N}\sum^N_{l=1}B^\eta_{ij}(X^{N,k}_{\eta,i}-X^{N,l}_{\eta,j})\Big)
-f'_{\gamma}\Big(\frac{1}{N}\sum^N_{l=1}B^\eta_{ij}(\bar X^{k}_{\eta,i}-\bar X^{l}_{\eta,j})\Big) \Big]\nonumber\\
&\qquad\qquad\qquad\qquad\qquad\qquad\cdot\frac{1}{N}\sum^N_{l=1}\nabla B^\eta_{ij}(\bar X^k_{\eta,i}-\bar X^l_{\eta,j})
\Big|^\kappa\,ds \Big)\nonumber\\
&+C\mathbb{E}\Big(N^{\alpha\kappa}\sum^n_{i=1}\max_{1\leq k\leq N}\int_0^{t\wedge\tau_\alpha}\Big|\sum^n_{j=1}
\Big[f'_\gamma\Big(\frac{1}{N}\sum^N_{l=1}B^\eta_{ij}(\bar X^k_{\eta,i}-\bar X^l_{\eta,j}) \Big)
-f'_\gamma(B^\eta_{ij}\ast u_{\eta,j}(\bar X^k_{\eta,i}) ) \Big]\nonumber\\
&\qquad\qquad\qquad\qquad\qquad\qquad\cdot\frac{1}{N}\sum^N_{l=1}\nabla B^\eta_{ij}(\bar X^k_{\eta,i}-\bar X^l_{\eta,j}) \Big|^\kappa\,ds \Big)\nonumber\\
&+C\mathbb{E}\Big(N^{\alpha\kappa}\sum^n_{i=1}\max_{1\leq k\leq N}\int_0^{t\wedge\tau_\alpha}
\Big|\sum^n_{j=1}f'_\gamma(B^\eta_{ij}\ast u_{\eta,j}(\bar X^k_{\eta,i}))\nonumber\\
&\qquad\qquad\qquad\qquad\qquad\qquad\cdot\Big[\frac{1}{N}\sum^N_{l=1}\nabla B^\eta_{ij}(\bar X^k_{\eta,i}-\bar X^l_{\eta,j})
-\nabla B^\eta_{ij}\ast u_{\eta,j}(\bar X^k_{\eta,i})  \Big]
 \Big|^\kappa\,ds \Big)\nonumber\\
=: & J^{21}+J^{22}+J^{23}+J^{24}.
\end{align}

The term $J^{21}$ can be divided into two terms:
\begin{align}\label{7}
J^{21}\leq &C\mathbb{E}\Big(N^{\alpha\kappa}\sum^n_{i=1}\max_{1\leq k\leq N}\int_0^{t\wedge\tau_\alpha}
\Big|\sum^n_{j=1}\Big[f'_\gamma\Big(\frac{1}{N}\sum^N_{l=N}B^\eta_{ij}(X^{N,k}_{\eta,i}-X^{N,l}_{\eta,j}) \Big)
-f'_\gamma\Big(\frac{1}{N}\sum^N_{l=N}B^\eta_{ij}(\bar X^{k}_{\eta,i}-\bar X^{l}_{\eta,j}) \Big) \Big]\nonumber\\
&\qquad\qquad\qquad\qquad\qquad\qquad \cdot\frac{1}{N}\sum^N_{l=1}[\nabla B^\eta_{ij}(X^{N,k}_{\eta,i}-X^{N,l}_{\eta,j})
-\nabla B^\eta_{ij}(\bar X^{k}_{\eta,i}-\bar X^{l}_{\eta,j}) ] \Big|^\kappa\,ds \Big)\nonumber\\
&+C\mathbb{E}\Big(N^{\alpha\kappa}\sum^n_{i=1}\max_{1\leq k\leq N}\int_0^{t\wedge\tau_\alpha}
\Big|\sum^n_{j=1}f'_\gamma\Big(\frac{1}{N}\sum^N_{l=N}B^\eta_{ij}(\bar X^{k}_{\eta,i}-\bar X^{l}_{\eta,j}) \Big) \nonumber\\
&\qquad\qquad\qquad\qquad\qquad\qquad \cdot\frac{1}{N}\sum^N_{l=1}[\nabla B^\eta_{ij}(X^{N,k}_{\eta,i}-X^{N,l}_{\eta,j})
-\nabla B^\eta_{ij}(\bar X^{k}_{\eta,i}-\bar X^{l}_{\eta,j}) ] \Big|^\kappa\,ds \Big)\nonumber\\
=:& J^{211}+J^{212}.
\end{align}
The term $J^{211}$ can be handled with
\begin{align}\label{8}
J^{211}&\leq C\|f''_\gamma\|^\kappa_{L^\infty(0,\sup_{1\leq i,j\leq n}\|B^\eta_{ij}\|_{L^\infty(\mathbb{R}^d)})}\sup_{1\leq i,j\leq n}\|\nabla B^\eta_{ij}\|^\kappa_{L^\infty(\mathbb{R}^d)}
\sup_{1\leq i,j\leq n}\|D^2B^\eta_{ij}\|^\kappa_{L^\infty(\mathbb{R}^d)}\nonumber\\
&\qquad\qquad\cdot\int_0^{t\wedge\tau_\alpha}\mathbb{E}
\Big(N^{\alpha\kappa}\sum^n_{i=1}\max_{1\leq k\leq N}|X^{N,k}_{\eta,i}-\bar X^k_{\eta,i}|^{2\kappa} \Big)\,ds\nonumber\\
&\leq C\|f''_\gamma\|^\kappa_{L^\infty(0,\sup_{1\leq i,j\leq n}\|B^\eta_{ij}\|_{L^\infty(\mathbb{R}^d)})}
\sup_{1\leq i,j\leq n}\eta^{-\kappa(\vartheta_{ij}+1)}
\sup_{1\leq i,j\leq n}\eta^{-\kappa(\vartheta_{ij}+2)}\nonumber\\
&\qquad\qquad\cdot N^{-\alpha\kappa}\int_0^{t\wedge\tau_\alpha}\mathbb{E}(S^\kappa_\alpha(s))\,ds\nonumber\\
&\leq C\sup_{1\leq i,j\leq n}N^{\kappa[\beta(2\vartheta_{ij}+4)-\alpha]}\int_0^{t}\mathbb{E}(S^\kappa_\alpha(s))\,ds,
\end{align}
where we have used the assumption $\eta=N^{-\beta}$ and
the fact that $\|f_\gamma\|_{C^3(\mathbb{R})}\leq \|f\|_{C^3([-\frac{1}{\gamma},\frac{1}{\gamma}])}\leq \frac{1}{\gamma^m}\leq N^\beta$.
With the help of $\|B^\eta_{ij}\ast u_{\eta,j}\|_{L^\infty((0,T)\times\mathbb{R}^d)}\leq C$,
we have\\
$\|f'_\gamma\|_{L^\infty(0,\sup_{1\leq i,j\leq n}\|B^\eta_{ij}\ast u_{\eta,j}\|_{L^\infty((0,T)\times\mathbb{R}^d)})}\leq C$.
Therefore,
\begin{align}\label{9}
J^{212}\leq& C\mathbb{E}\Big(N^{\alpha\kappa}\sum^n_{i=1}\max_{1\leq k\leq N}\int_0^{t\wedge\tau_\alpha}
\Big|\sum^n_{j=1}\Big[f'_\gamma\Big(\frac{1}{N}\sum^N_{l=1}B^\eta_{ij}(\bar X^k_{\eta,i}-\bar X^l_{\eta,j}) \Big)
-f'_\gamma(B^\eta_{ij}\ast u_{\eta,j}(\bar X^k_{\eta,i}) ) \Big]\nonumber\\
&\qquad\qquad\qquad\qquad\qquad\qquad \cdot\frac{1}{N}\sum^N_{l=1}[\nabla B^\eta_{ij}(X^{N,k}_{\eta,i}-X^{N,l}_{\eta,j})
-\nabla B^\eta_{ij}(\bar X^{k}_{\eta,i}-\bar X^{l}_{\eta,j}) ] \Big|^\kappa\,ds \Big)\nonumber\\
&+C\mathbb{E}\Big(N^{\alpha\kappa}\sum^n_{i=1}\max_{1\leq k\leq N}\int_0^{t\wedge\tau_\alpha}
 \Big|\sum^n_{j=1}f'_\gamma(B^\eta_{ij}\ast u_{\eta,j}(\bar X^k_{\eta,i}) )\nonumber\\
 &\qquad\qquad\qquad\qquad\qquad\qquad \cdot\frac{1}{N}\sum^N_{l=1}[\nabla B^\eta_{ij}(X^{N,k}_{\eta,i}-X^{N,l}_{\eta,j})
-\nabla B^\eta_{ij}(\bar X^{k}_{\eta,i}-\bar X^{l}_{\eta,j}) ] \Big|^\kappa\,ds \Big)\nonumber\\
\leq &C\|f''_\gamma\|^\kappa_{L^\infty(0,\sup_{1\leq i,j\leq n}\|B^\eta_{ij}\|_{L^\infty(\mathbb{R}^d)})}
\sup_{1\leq i,j\leq n}\eta^{-\kappa(\vartheta_{ij}+2)}\nonumber\\
&\qquad \cdot \mathbb{E}\Big(\int_0^{t\wedge\tau_\alpha}\sum^n_{i,j=1}\max_{1\leq k\leq N}
\Big|\frac{1}{N}\sum^N_{l=1}B^\eta_{ij}(\bar X^k_{\eta,i}-\bar X^l_{\eta,j})-B^\eta_{ij}\ast u_{\eta,j}(\bar X^k_{\eta,i}) \Big|^\kappa
S^\kappa_\alpha(s)\,ds \Big)\nonumber\\
&+C \mathbb{E}\Big(N^{\alpha\kappa}\sum^n_{i=1}\max_{1\leq k\leq N}\int_0^{t\wedge\tau_\alpha}
\sum^n_{j=1}\Big|\frac{1}{N}\sum^N_{l=1}D^2B^\eta_{ij}(\bar X^k_{\eta,i}-\bar X^l_{\eta,j})
(X^{N,k}_{\eta,i}-\bar X^k_{\eta,i}) \Big|^\kappa\,ds  \Big)\nonumber\\
&+C \mathbb{E}\Big(N^{\alpha\kappa}\sum^n_{i=1}\max_{1\leq k\leq N}\int_0^{t\wedge\tau_\alpha}
\sum^n_{j=1}\Big|\frac{1}{N}\sum^N_{l=1}D^2B^\eta_{ij}(\bar X^k_{\eta,i}-\bar X^l_{\eta,j})
(X^{N,l}_{\eta,j}-\bar X^l_{\eta,j}) \Big|^\kappa\,ds  \Big)\nonumber\\
&+C\sup_{1\leq i,j\leq n}\|D^3B^\eta_{ij}\|^\kappa_{L^\infty(\mathbb{R}^d)} \mathbb{E}\Big(\int_0^{t\wedge \tau_\alpha}N^{\alpha\kappa}\sum^n_{i=1}
\max_{1\leq k\leq N} |X^{N,k}_{\eta,i}-\bar X^k_{\eta,i}|^{2\kappa}\,ds \Big)\nonumber\\
=:&J^{2121}+J^{2122}+J^{2123}
+C\sup_{1\leq i,j\leq n}N^{\kappa[\beta(\vartheta_{ij}+3)-\alpha]}\int_0^{t} \mathbb{E}(S^\kappa_\alpha(s))\,ds.
\end{align}
For $J^{2121}$, we split the domain $\Omega=A^{N,n}_{\theta_1,B^\eta}\cup (A^{N,n}_{\theta_1,B^\eta})^c$ and obtain
\begin{align}\label{10}
J^{2121}\leq& C\|f''_\gamma\|^\kappa_{L^\infty(0,\sup_{1\leq i,j\leq n}\|B^\eta_{ij}\|_{L^\infty(\mathbb{R}^d)})}
\sup_{1\leq i,j\leq n}\eta^{-\kappa(\vartheta_{ij}+2)}\nonumber\\
&\qquad\cdot \mathbb{E}\Big(\int_0^{t\wedge\tau_\alpha}\sum^n_{i,j=1}\max_{1\leq k\leq N}
\Big|\frac{1}{N}\sum^N_{l=1}B^\eta_{ij}(\bar X^k_{\eta,i}-\bar X^l_{\eta,j})-B^\eta_{ij}\ast u_{\eta,j}(\bar X^k_{\eta,i}) \Big|^\kappa
\mathbb{I}_{A^{N,n}_{\theta_1,B^\eta}} S^\kappa_\alpha(s)\,ds \Big)\nonumber\\
&+ C\|f''_\gamma\|^\kappa_{L^\infty(0,\sup_{1\leq i,j\leq n}\|B^\eta_{ij}\|_{L^\infty(\mathbb{R}^d)})}
\sup_{1\leq i,j\leq n}\eta^{-\kappa(\vartheta_{ij}+2)}\nonumber\\
&\qquad\cdot \mathbb{E}\Big(\int_0^{t\wedge\tau_\alpha}\sum^n_{i,j=1}\max_{1\leq k\leq N}
\Big|\frac{1}{N}\sum^N_{l=1}B^\eta_{ij}(\bar X^k_{\eta,i}-\bar X^l_{\eta,j})-B^\eta_{ij}\ast u_{\eta,j}(\bar X^k_{\eta,i}) \Big|^\kappa
\mathbb{I}_{(A^{N,n}_{\theta_1,B^\eta})^c} S^\kappa_\alpha(s)\,ds \Big)\nonumber\\
\leq& C\sup_{1\leq i,j\leq n}N^{\kappa\beta(\vartheta_{ij}+3)}\big(\sup_{1\leq i,j\leq n}\|B^\eta_{ij}\|^\kappa_{L^\infty(\mathbb{R}^d)}
+\sup_{1\leq i,j\leq n}\|B^\eta_{ij}\ast u_{\eta,j}\|^\kappa_{L^\infty((0,T)\times\mathbb{R}^d)}\big)
\int_0^t \mathbb{P}(A^{N,n}_{\theta_1,B^\eta})\,ds\nonumber\\
&+ C\sup_{1\leq i,j\leq n}N^{\kappa[\beta(\vartheta_{ij}+3)-\theta_1]}\int_0^{t\wedge\tau_\alpha} \mathbb{E}(S^\kappa_\alpha(s))\,ds\nonumber\\
\leq & C\sup_{1\leq i,j\leq n}N^{\kappa\beta(2\vartheta_{ij}+3)}\int_0^t \mathbb{P}(A^{N,n}_{\theta_1,B^\eta})\,ds
+C\sup_{1\leq i,j\leq n}N^{\kappa[\beta(\vartheta_{ij}+3)-\theta_1]}\int_0^t \mathbb{E}(S_\alpha^\kappa(s))\,ds,
\end{align}
where we have used $\|f_\gamma\|_{C^3(\mathbb{R})}\leq \frac{1}{\gamma^m}\leq N^\beta$
and $\|B^\eta_{ij}\|_{L^\infty(\mathbb{R}^d)}\leq \frac{C}{\eta^{\vartheta_{ij}}}\leq CN^{\beta\vartheta_{ij}}$,
more details can be found in \cite{CGHL}.
For $J^{2123}$, we split again the domain $\Omega={A^{N,n}_{0,|D^2B^\eta|}}\cup ({A^{N,n}_{0,|D^2B^\eta|}})^c$ and obtain
\begin{align}\label{11}
J^{2123}\leq& C \mathbb{E}\Big(N^{\alpha\kappa}\sum^n_{i=1}\max_{1\leq k\leq N}\int_0^{t\wedge\tau_\alpha}
\sum^n_{j=1}\Big|\frac{1}{N}\sum^N_{l=1}|D^2B^\eta_{ij}|(\bar X^k_{\eta,i}-\bar X^l_{\eta,j}) \Big|^\kappa
\sup_l|X^{N,l}_{\eta,j}-\bar X^l_{\eta,j}|^\kappa\,ds \Big)\nonumber\\
\leq&C \mathbb{E}\Big(\int_0^{t\wedge\tau_\alpha}
\sum^n_{i,j=1}\max_{1\leq k\leq N}\Big|\frac{1}{N}\sum^N_{l=1}|D^2B^\eta_{ij}|(\bar X^k_{\eta,i}-\bar X^l_{\eta,j})
-|D^2B^\eta_{ij}|\ast u_{\eta,j}(\bar X^k_{\eta,i}) \Big|^\kappa
\mathbb{I}_{A^{N,n}_{0,|D^2B^\eta|}}S^\kappa_\alpha(s) \,ds  \Big)\nonumber\\
&+C \mathbb{E}\Big(\int_0^{t\wedge\tau_\alpha}
\sum^n_{i,j=1}\max_{1\leq k\leq N}\Big|\frac{1}{N}\sum^N_{l=1}|D^2B^\eta_{ij}|(\bar X^k_{\eta,i}-\bar X^l_{\eta,j})
-|D^2B^\eta_{ij}|\ast u_{\eta,j}(\bar X^k_{\eta,i}) \Big|^\kappa
\mathbb{I}_{(A^{N,n}_{0,|D^2B^\eta|})^c}S^\kappa_\alpha(s) \,ds  \Big)\nonumber\\
&+C \sup_{1\leq i,j\leq n}\big\| |D^2B^\eta_{ij}|\ast u_{\eta,j}\big\|^\kappa_{L^\infty((0,T)\times\mathbb{R}^d)}
\int_0^{t\wedge\tau_\alpha} \mathbb{E}(S^\kappa_\alpha(s))\,ds\nonumber\\
\leq & C\sup_{1\leq i,j\leq n}N^{\kappa\beta(\vartheta_{ij}+2)}\int_0^t \mathbb{P}(A^{N,n}_{0,|D^2B^\eta|})\,ds
-C\ln\eta\int_0^{t} \mathbb{E}(S^\kappa_\alpha(s))\,ds,
\end{align}
where we have used the results that
\begin{align*}
\big\||D^2B^\eta_{ij}|\ast u_{\eta,j}\big\|_{L^\infty((0,T)\times\mathbb{R}^d)}\leq
\left\{
       \begin{array}{lcr}
       C  \;\;\;{\rm{if}}\;\;\;0<\vartheta_{ij}< d-2,\\
       -C\ln\eta\;\;\;{\rm{if}}\;\;\;\vartheta_{ij}=d-2,
       \end{array}
\right.
\end{align*}
more details can be found in \cite{CGHL}.
Similarly, we can derive that
\begin{align}\label{12}
J^{2122}\leq  C\sup_{1\leq i,j\leq n}N^{\kappa\beta(\vartheta_{ij}+2)}\int_0^t \mathbb{P}(A^{N,n}_{0,D^2B^\eta})\,ds
+C\int_0^{t} \mathbb{E}(S^\kappa_\alpha(s))\,ds,
\end{align}
where we have used $\|D^2 B^\eta_{ij}\ast u_{\eta,j}\|_{L^\infty((0,T)\times\mathbb{R}^d)}\leq C$.
Combining \eqref{7}-\eqref{12}, we have
\begin{align}\label{14}
J^{21}\leq&
C\sup_{1\leq i,j\leq n}\big(1+\ln N+N^{\kappa[\beta(2\vartheta_{ij}+4)-\alpha]}+N^{\kappa[\beta(\vartheta_{ij}+3)-\alpha]}
+N^{\kappa[\beta(\vartheta_{ij}+3)-\theta_1]}\big)\int_0^{t}\mathbb{E}(S^\kappa_\alpha(s))\,ds\nonumber\\
&+C\sup_{1\leq i,j\leq n}N^{\kappa\beta(2\vartheta_{ij}+3)}\int_0^t \mathbb{P}(A^{N,n}_{\theta_1,B^\eta})\,ds
+C\sup_{1\leq i,j\leq n}N^{\kappa\beta(\vartheta_{ij}+2)}\int_0^t \mathbb{P}(A^{N,n}_{0,D^2B^\eta})\,ds\nonumber\\
&+C\sup_{1\leq i,j\leq n}N^{\kappa\beta(\vartheta_{ij}+2)}\int_0^t \mathbb{P}(A^{N,n}_{0,|D^2B^\eta|})\,ds.
\end{align}

For $J^{23}$, we split again the domain $\Omega={A^{N,n}_{\theta_2,B^\eta}}\cup (A^{N,n}_{\theta_2,B^\eta})^c$ and obtain
\begin{align}\label{23}
J^{23}\leq &C\sup_{1\leq i,j\leq n}\|\nabla B^\eta_{ij}\|^\kappa_{L^\infty(\mathbb{R}^d)}
\|f''_\gamma\|^\kappa_{L^\infty(0,\sup_{1\leq i,j\leq n}\|B^\eta_{ij}\|_{L^\infty(\mathbb{R}^d)})}\nonumber\\
&\qquad \cdot \mathbb{E}\Big(\int_0^{t\wedge\tau_\alpha}N^{\alpha\kappa}\sum^n_{i,j=1}\max_{1\leq k\leq N}
\Big|\frac{1}{N}\sum^N_{l=1}B^\eta_{ij}(\bar X^k_{\eta,i}-\bar X^l_{\eta,j})-B^\eta_{ij}\ast u_{\eta,j}(\bar X^k_{\eta,i}) \Big|^\kappa
\mathbb{I}_{A^{N,n}_{\theta_2,B^\eta}} \,ds \Big)\nonumber\\
&+C\sup_{1\leq i,j\leq n}\|\nabla B^\eta_{ij}\|^\kappa_{L^\infty(\mathbb{R}^d)}
\|f''_\gamma\|^\kappa_{L^\infty(0,\sup_{1\leq i,j\leq n}\|B^\eta_{ij}\|_{L^\infty(\mathbb{R}^d)})}\nonumber\\
&\qquad \cdot \mathbb{E}\Big(\int_0^{t\wedge\tau_\alpha}N^{\alpha\kappa}\sum^n_{i,j=1}\max_{1\leq k\leq N}
\Big|\frac{1}{N}\sum^N_{l=1}B^\eta_{ij}(\bar X^k_{\eta,i}-\bar X^l_{\eta,j})-B^\eta_{ij}\ast u_{\eta,j}(\bar X^k_{\eta,i}) \Big|^\kappa
\mathbb{I}_{(A^{N,n}_{\theta_2,B^\eta})^c} \,ds \Big)\nonumber\\
\leq & C\sup_{1\leq i,j\leq n}N^{\kappa[\beta(2\vartheta_{ij}+2)+\alpha]}\int_0^t \mathbb{P}(A^{N,n}_{\theta_2,B^\eta})\,ds
+C\sup_{1\leq i,j\leq n}N^{\kappa[\beta(\vartheta_{ij}+2)+\alpha-\theta_2]}.
\end{align}

For $J^{24}$, we split again the domain $\Omega={A^{N,n}_{\theta_3,\nabla B^\eta}}\cup (A^{N,n}_{\theta_3,\nabla B^\eta})^c$ and obtain
\begin{align}\label{24}
J^{24}\leq &CN^{\alpha\kappa} \mathbb{E}\Big(\int_0^{t\wedge\tau_\alpha}\sum^n_{i,j=1}\max_{1\leq k\leq N}
\Big|\frac{1}{N}\sum^N_{l=1}\nabla B^\eta_{ij}(\bar X^k_{\eta,i}-\bar X^l_{\eta,j})-\nabla B^\eta_{ij}\ast u_{\eta,j}(\bar X^k_{\eta,i}) \Big|^\kappa
\mathbb{I}_{A^{N,n}_{\theta_3,\nabla B^\eta}}\,ds \Big)\nonumber\\
&+CN^{\alpha\kappa} \mathbb{E}\Big(\int_0^{t\wedge\tau_\alpha}\sum^n_{i,j=1}\max_{1\leq k\leq N}
\Big|\frac{1}{N}\sum^N_{l=1}\nabla B^\eta_{ij}(\bar X^k_{\eta,i}-\bar X^l_{\eta,j})-\nabla B^\eta_{ij}\ast u_{\eta,j}(\bar X^k_{\eta,i}) \Big|^\kappa
\mathbb{I}_{(A^{N,n}_{\theta_3,\nabla B^\eta})^c}\,ds \Big)\nonumber\\
\leq &C\sup_{1\leq i,j\leq n}N^{\kappa[\beta(\vartheta_{ij}+1)+\alpha]}\int_0^t \mathbb{P}(A^{N,n}_{\theta_3,\nabla B^\eta})\,ds+CN^{\kappa(\alpha-\theta_3)}.
\end{align}

The term $J^{22}$ can be divided into two terms.
\begin{align}\label{15}
J^{22}\leq &
C\mathbb{E}\Big(N^{\alpha\kappa}\sum^n_{i=1}\max_{1\leq k\leq N}\int_0^{t\wedge \tau_\alpha}
\sum^n_{j=1}\Big|f''_\gamma\Big(\frac{1}{N}\sum^N_{l=1}B^\eta_{ij}(\bar X^k_{\eta,i}-\bar X^l_{\eta,j}) \Big)\nonumber\\
&\qquad\qquad\cdot \frac{1}{N}\sum^N_{l=1}[B^\eta_{ij}(X^{N,k}_{\eta,i}-X^{N,l}_{\eta,j})-B^\eta_{ij}(\bar X^k_{\eta,i}-\bar X^l_{\eta,j})]  \Big|^\kappa
  \Big|\frac{1}{N}\sum^N_{l=1}\nabla B^\eta_{ij}(\bar X^k_{\eta,i}-\bar X^l_{\eta,j}) \Big|^\kappa\,ds \Big)\nonumber\\
&+C\mathbb{E}\Big(N^{\alpha\kappa}\sum^n_{i=1}\max_{1\leq k\leq N}\int_0^{t\wedge \tau_\alpha}\|f'''_\gamma\|^\kappa_{L^\infty(0,\sup_{1\leq i,j\leq n}\|B^\eta_{ij}\|_{L^\infty(\mathbb{R}^d)})}\nonumber\\
&\qquad \cdot\sum^n_{j=1}\Big|\frac{1}{N}\sum^N_{l=1}[B^\eta_{ij}(X^{N,k}_{\eta,i}-X^{N,l}_{\eta,j})-B^\eta_{ij}(\bar X^k_{\eta,i}-\bar X^l_{\eta,j})]  \Big|^{2\kappa}
 \Big|\frac{1}{N}\sum^N_{l=1}\nabla B^\eta_{ij}(\bar X^k_{\eta,i}-\bar X^l_{\eta,j}) \Big|^\kappa\,ds \Big)\nonumber\\
=: &J^{221}+J^{222}.
\end{align}
For $J^{222}$,
\begin{align}\label{16}
J^{222}\leq&C\|f'''_\gamma\|^\kappa_{L^\infty(0,\sup_{1\leq i,j\leq n}\|B^\eta_{ij}\|_{L^\infty(\mathbb{R}^d)})} \sup_{1\leq i,j\leq n}\|\nabla B^\eta_{ij}\|^{3\kappa}_{L^\infty(\mathbb{R}^d)} \mathbb{E}\Big(\int_0^{t\wedge \tau_\alpha}N^{\alpha\kappa}\sum^n_{i=1}\max_{1\leq k\leq N} |X^{N,k}_{\eta,i}-\bar X^k_{\eta,i}|^{2\kappa}\,ds \Big)\nonumber\\
\leq & C\sup_{1\leq i,j\leq n}N^{\kappa[\beta(3\vartheta_{ij}+4)-\alpha]}\int_0^t \mathbb{E}(S^\kappa_\alpha(s))\,ds.
\end{align}
Now we focus on $J^{221}$.
\begin{align}\label{17}
J^{221}\leq& C\mathbb{E}\Big(N^{\alpha\kappa}\sum^n_{i=1}\max_{1\leq k\leq N}\int_0^{t\wedge\tau_\alpha}
\sum^n_{j=1}\Big|f''_\gamma\Big(\frac{1}{N}\sum^N_{l=1}B^\eta_{ij}(\bar X^k_{\eta,i}-\bar X^l_{\eta,j}) \Big) \Big|^\kappa\nonumber\\
&\qquad \cdot\Big|\frac{1}{N}\sum^N_{l=1}|\nabla B^\eta_{ij}(\bar X^k_{\eta,i}-\bar X^l_{\eta,j})|\cdot|X^{N,l}_{\eta,j}-\bar X^l_{\eta,j}|\Big|^\kappa
 \Big|\frac{1}{N}\sum^N_{l=1}\nabla B^\eta_{ij}(\bar X^k_{\eta,i}-\bar X^l_{\eta,j}) \Big|^\kappa \,ds\Big)\nonumber\\
&+C\|f''_\gamma\|^\kappa_{L^\infty(0,\sup_{1\leq i,j\leq n}\|B^\eta_{ij}\|_{L^\infty(\mathbb{R}^d)})}\sup_{1\leq i,j\leq n}\|\nabla B^\eta_{ij}\|^\kappa_{L^\infty(\mathbb{R}^d)}
\sup_{1\leq i,j\leq n}\|D^2B^\eta_{ij}\|^\kappa_{L^\infty(\mathbb{R}^d)}\nonumber\\
&\qquad \cdot \mathbb{E}\Big(\int_0^{t\wedge\tau_\alpha}N^{\alpha\kappa}\sum^n_{i=1}\max_{1\leq k\leq N}|X^{N,k}_{\eta,i}-\bar X^k_{\eta,i}|^{2\kappa} \,ds \Big)\nonumber\\
\leq &C\mathbb{E}\Big(N^{\alpha\kappa}\sum^n_{i=1}\max_{1\leq k\leq N}\int_0^{t\wedge\tau_\alpha}
\sum^n_{j=1}\Big|f''_\gamma\Big(\frac{1}{N}\sum^N_{l=1}B^\eta_{ij}(\bar X^k_{\eta,i}-\bar X^l_{\eta,j}) \Big)
-f''_\gamma(B^\eta_{ij}\ast u_{\eta,i}(\bar X^k_{\eta,i})) \Big|^\kappa\nonumber\\
&\qquad\cdot \|\nabla B^\eta_{ij}\|^{2\kappa}_{L^\infty(\mathbb{R}^d)}
\max_{1\leq l\leq N}|X^{N,l}_{\eta,j}-\bar X^l_{\eta,j}|^\kappa\,ds \Big)\nonumber\\
&+C\mathbb{E}\Big(N^{\alpha\kappa}\sum^n_{i=1}\max_{1\leq k\leq N}\int_0^{t\wedge\tau_\alpha}
\sum^n_{j=1}|f''_\gamma(B^\eta_{ij}\ast u_{\eta,j}(\bar X^k_{\eta,i}))|^\kappa\nonumber\\
&\qquad \cdot\Big|\frac{1}{N}\sum^N_{l=1}
|\nabla B^\eta_{ij}(\bar X^k_{\eta,i}-\bar X^l_{\eta,j})|-|\nabla B^\eta_{ij}|\ast u_{\eta,j}(\bar X^k_{\eta,i}) \Big|^\kappa
\max_{1\leq l\leq N}|X^{N,l}_{\eta,j}-\bar X^l_{\eta,j}|^\kappa\|\nabla B^\eta_{ij}\|^\kappa_{L^\infty(\mathbb{R}^d)}\,ds \Big)\nonumber\\
&+C\mathbb{E}\Big(N^{\alpha\kappa}\sum^n_{i=1}\max_{1\leq k\leq N}\int_0^{t\wedge\tau_\alpha}
\sum^n_{j=1}|f''_\gamma(B^\eta_{ij}\ast u_{\eta,j}(\bar X^k_{\eta,i}))|^\kappa
\big| |\nabla B^\eta_{ij}|\ast u_{\eta,j}(\bar X^k_{\eta,i})\big|^\kappa\nonumber\\
&\qquad\quad \cdot\max_{1\leq l\leq N}|X^{N,l}_{\eta,j}-\bar X^l_{\eta,j}|^\kappa
 \Big|\frac{1}{N}\sum^N_{l=1}\nabla B^\eta_{ij}(\bar X^k_{\eta,i}-\bar X^l_{\eta,j})
-\nabla B^\eta_{ij}\ast u_{\eta,j}(\bar X^k_{\eta,i}) \Big|^\kappa \,ds \Big)\nonumber\\
&+C\mathbb{E}\Big(N^{\alpha\kappa}\sum^n_{i=1}\max_{1\leq k\leq N}\int_0^{t\wedge\tau_\alpha}
\sum^n_{j=1}|f''_\gamma(B^\eta_{ij}\ast u_{\eta,j}(\bar X^k_{\eta,i}))|^\kappa
\big| |\nabla B^\eta_{ij}|\ast u_{\eta,j}(\bar X^k_{\eta,i})\big|^\kappa\nonumber\\
&\qquad\quad \cdot\max_{1\leq l\leq N}|X^{N,l}_{\eta,j}-\bar X^l_{\eta,j}|^\kappa
 \Big|\nabla B^\eta_{ij}\ast u_{\eta,j}(\bar X^k_{\eta,i}) \Big|^\kappa \,ds \Big)\nonumber\\
&+C\|f''_\gamma\|^\kappa_{L^\infty(0,\sup_{1\leq i,j\leq n}\|B^\eta_{ij}\|_{L^\infty(\mathbb{R}^d)})}
N^{\kappa[\beta(\vartheta_{ij}+1)+\beta(\vartheta_{ij}+2)-\alpha]}\int_0^t \mathbb{E}(S^\kappa_\alpha(s))\,ds\nonumber\\
=:&J^{2211}+J^{2212}+J^{2213}+C\int_0^t \mathbb{E}(S^\kappa_\alpha(s))\,ds
+C\sup_{1\leq i,j\leq n}N^{\kappa[\beta(2\vartheta_{ij}+4)-\alpha]}\int_0^t \mathbb{E}(S^\kappa_\alpha(s))\,ds,
\end{align}
where we have used
$\|\nabla B^\eta_{ij}\ast u_{\eta,j}\|_{L^\infty((0,T)\times\mathbb{R}^d)}
+\big\||\nabla B^\eta_{ij}|\ast u_{\eta,j} \big\|_{L^\infty((0,T)\times\mathbb{R}^d)}\leq C$.
For $J^{2211}$, we split again the domain $\Omega={A^{N,n}_{\theta_4,B^\eta}}\cup (A^{N,n}_{\theta_4,B^\eta})^c$ and obtain
\begin{align}\label{18}
J^{2211}\leq &C\sup_{1\leq i,j\leq n}N^{\kappa\beta(2\vartheta_{ij}+3)} \mathbb{E}\Big(\int_0^{t\wedge\tau_\alpha}\sum^n_{i,j=1}\max_{1\leq k\leq N}
\Big|\frac{1}{N}\sum^N_{l=1}B^\eta_{ij}(\bar X^k_{\eta,i}-\bar X^l_{\eta,j})-B^\eta_{ij}\ast u_{\eta,j}(\bar X^k_{\eta,i}) \Big|^\kappa
\mathbb{I}_{A^{N,n}_{\theta_4,B^\eta}}S^\kappa_\alpha(s)\,ds \Big)\nonumber\\
&+C\sup_{1\leq i,j\leq n}N^{\kappa\beta(2\vartheta_{ij}+3)} \mathbb{E}\Big(\int_0^{t\wedge\tau_\alpha}\sum^n_{i,j=1}\max_{1\leq k\leq N}
\Big|\frac{1}{N}\sum^N_{l=1}B^\eta_{ij}(\bar X^k_{\eta,i}-\bar X^l_{\eta,j})-B^\eta_{ij}\ast u_{\eta,j}(\bar X^k_{\eta,i}) \Big|^\kappa
\mathbb{I}_{(A^{N,n}_{\theta_4,B^\eta})^c}S^\kappa_\alpha(s)\,ds \Big)\nonumber\\
\leq & C\sup_{1\leq i,j\leq n}N^{\kappa\beta(3\vartheta_{ij}+3)}\int_0^t \mathbb{P}(A^{N,n}_{\theta_4,B^\eta})\,ds
+C\sup_{1\leq i,j\leq n}N^{\kappa[\beta(2\vartheta_{ij}+3)-\theta_4]}\int_0^t \mathbb{E}(S^\kappa_\alpha(s))\,ds.
\end{align}
For $J^{2212}$, we split again the domain $\Omega={A^{N,n}_{\theta_5,|\nabla B^\eta|}}\cup (A^{N,n}_{\theta_5,|\nabla B^\eta|})^c$ and obtain
\begin{align}\label{19}
&J^{2212}\nonumber\\
\leq& C\sup_{1\leq i,j\leq n}N^{\kappa\beta(\vartheta_{ij}+1)} \mathbb{E}\Big(\int_0^{t\wedge\tau_\alpha}\sum^n_{i,j=1}\max_{1\leq k\leq N}
\Big|\frac{1}{N}\sum^N_{l=1}|\nabla B^\eta_{ij}(\bar X^k_{\eta,i}-\bar X^l_{\eta,j})|
-|\nabla B^\eta_{ij}|\ast u_{\eta,j}(\bar X^k_{\eta,i}) \Big|^\kappa\mathbb{I}_{A^{N,n}_{\theta_5,|\nabla B^\eta|}}S^\kappa_\alpha(s)\,ds \Big)\nonumber\\
&+C\sup_{1\leq i,j\leq n}N^{\kappa\beta(\vartheta_{ij}+1)} \mathbb{E}\Big(\int_0^{t\wedge\tau_\alpha}\sum^n_{i,j=1}\max_{1\leq k\leq N}
\Big|\frac{1}{N}\sum^N_{l=1}|\nabla B^\eta_{ij}(\bar X^k_{\eta,i}-\bar X^l_{\eta,j})|
-|\nabla B^\eta_{ij}|\ast u_{\eta,j}(\bar X^k_{\eta,i}) \Big|^\kappa\mathbb{I}_{(A^{N,n}_{\theta_5,|\nabla B^\eta|})^c}S^\kappa_\alpha(s)\,ds \Big)\nonumber\\
\leq& C\sup_{1\leq i,j\leq n}N^{2\kappa\beta(\vartheta_{ij}+1)}\int^t_0 \mathbb{P}(A^{N,n}_{\theta_5,|\nabla B^\eta|})\,ds
+C\sup_{1\leq i,j\leq n}N^{\kappa[\beta(\vartheta_{ij}+1)-\theta_5]}\int_0^t \mathbb{E}(S^\kappa_\alpha(s))\,ds.
\end{align}
For $J^{2213}$, we take $\theta_6=0$ and obtain
\begin{align}\label{20}
J^{2213}\leq &C\mathbb{E}\Big(\int_0^{t\wedge\tau_\alpha}
\sum^n_{i,j=1}\max_{1\leq k\leq N} \Big|\frac{1}{N}\sum^N_{l=1}\nabla B^\eta_{ij}(\bar X^k_{\eta,i}-\bar X^l_{\eta,j})
-\nabla B^\eta_{ij}\ast u_{\eta,i}(\bar X^k_{\eta,i}) \Big|^\kappa\mathbb{I}_{A^{N,n}_{0,\nabla B^\eta}}S^\kappa_\alpha(s)\,ds \Big)\nonumber\\
&+C\mathbb{E}\Big(\int_0^{t\wedge\tau_\alpha}
\sum^n_{i,j=1}\max_{1\leq k\leq N} \Big|\frac{1}{N}\sum^N_{l=1}\nabla B^\eta_{ij}(\bar X^k_{\eta,i}-\bar X^l_{\eta,j})
-\nabla B^\eta_{ij}\ast u_{\eta,i}(\bar X^k_{\eta,i}) \Big|^\kappa\mathbb{I}_{(A^{N,n}_{0,\nabla B^\eta})^c}S^\kappa_\alpha(s)\,ds \Big)\nonumber\\
\leq & C\sup_{1\leq i,j\leq n}N^{\kappa\beta(\vartheta_{ij}+1)}\int_0^t \mathbb{P}(A^{N,n}_{0,\nabla B^\eta})\,ds+ C\int_0^t \mathbb{E}(S^\kappa_\alpha(s))\,ds.
\end{align}
From \eqref{15}-\eqref{20}, it holds
\begin{align}\label{22}
J^{22}\leq &
C\sup_{1\leq i,j\leq n}\big(1+N^{\kappa[\beta(\vartheta_{ij}+1)-\theta_5]}+N^{\kappa[\beta(2\vartheta_{ij}+3)-\theta_4]}+N^{\kappa[\beta(3\vartheta_{ij}+4)-\alpha]}\big)
\int_0^t \mathbb{E}(S^\kappa_\alpha(s))\,ds \nonumber\\
&+C\sup_{1\leq i,j\leq n}N^{\kappa\beta(3\vartheta_{ij}+3)}\int_0^t \mathbb{P}(A^{N,n}_{\theta_4,B^\eta})\,ds
+C\sup_{1\leq i,j\leq n}N^{\kappa\beta(2\vartheta_{ij}+2)}\int_0^t \mathbb{P}(A^{N,n}_{\theta_5,|\nabla B^\eta|})\,ds\nonumber\\
&+C\sup_{1\leq i,j\leq n}N^{\kappa\beta(\vartheta_{ij}+1)}\int_0^t \mathbb{P}(A^{N,n}_{0,\nabla B^\eta})\,ds.
\end{align}
Plugging \eqref{14}-\eqref{24} and \eqref{22} into \eqref{6}, and combining \eqref{5} and Lemma \ref{large}, we get
\begin{align}\label{25}
E(S^\kappa_\alpha(t))\leq &C\sup_{1\leq i,j\leq n}\big(1+\ln N+N^{\kappa[\beta(2\vartheta_{ij}+4)-\alpha]}
+N^{\kappa[\beta(\vartheta_{ij}+3)-\alpha]}
+N^{\kappa[\beta(\vartheta_{ij}+3)-\theta_1]}\nonumber\\
&\qquad+N^{\kappa[\beta(\vartheta_{ij}+1)-\theta_5]}+N^{\kappa[\beta(2\vartheta_{ij}+3)-\theta_4]}
+N^{\kappa[\beta(3\vartheta_{ij}+4)-\alpha]}\big)\int_0^t \mathbb{E}(S^\kappa_\alpha(s))\,ds\nonumber\\
&+C\sup_{1\leq i,j\leq n}N^{\kappa\beta(2\vartheta_{ij}+3)+2m_1(\theta_1-\frac{1}{2}+\beta \vartheta_{ij})+1}
+C\sup_{1\leq i,j\leq n}N^{\kappa\beta(\vartheta_{ij}+2)+2m_7[-\frac{1}{2}+\beta(\vartheta_{ij}+2)]+1}\nonumber\\
&+C\sup_{1\leq i,j\leq n}N^{\kappa\beta(\vartheta_{ij}+2)+2m_8[-\frac{1}{2}+\beta(\vartheta_{ij}+2)]+1}
+C\sup_{1\leq i,j\leq n}N^{\kappa[\beta(2\vartheta_{ij}+2)+\alpha]+2m_2(\theta_2-\frac{1}{2}+\beta \vartheta_{ij})+1}\nonumber\\
&+C\sup_{1\leq i,j\leq n}N^{\kappa[\beta(\vartheta_{ij}+1)+\alpha]+2m_3[\theta_3-\frac{1}{2}+\beta (\vartheta_{ij}+1)]+1}
+C\sup_{1\leq i,j\leq n}N^{\kappa\beta(3\vartheta_{ij}+3)+2m_4(\theta_4-\frac{1}{2}+\beta \vartheta_{ij})+1}\nonumber\\
&+C\sup_{1\leq i,j\leq n}N^{\kappa\beta(2\vartheta_{ij}+2)+2m_5[\theta_5-\frac{1}{2}+\beta (\vartheta_{ij}+1)]+1}
+C\sup_{1\leq i,j\leq n}N^{\kappa\beta(\vartheta_{ij}+1)+2m_6[-\frac{1}{2}+\beta (\vartheta_{ij}+1)]+1}\nonumber\\
&+C\sup_{1\leq i,j\leq n}N^{\kappa[\beta(\vartheta_{ij}+2)+\alpha-\theta_2]}+CN^{\kappa(\alpha-\theta_3)}.
\end{align}
For a given $\beta$, we choose $\alpha$ such that
\begin{align}\label{26}
\sup_{1\leq i,j\leq n}\beta(3\vartheta_{ij}+4)\leq \alpha.
\end{align}

To bound all the terms above, we need the following restrictions for $\theta_1$ - $\theta_5$.
\begin{align}
&\sup_{1\leq i,j\leq n}\beta(\vartheta_{ij}+3)\leq \theta_1<\frac{1}{2}-\beta\sup_{1\leq i,j\leq n} \vartheta_{ij},\label{27}\\
&0<\frac{1}{2}-\beta(\sup_{1\leq i,j\leq n}\vartheta_{ij}+2),\label{theta1}\\
&\sup_{1\leq i,j\leq n}\beta(\vartheta_{ij}+2)+\alpha<\theta_2<\frac{1}{2}-\beta \sup_{1\leq i,j\leq n}\vartheta_{ij},\label{28}\\
&\alpha<\theta_3<\frac{1}{2}-\beta(\sup_{1\leq i,j\leq n}\vartheta_{ij}+1),\label{29}\\
&\sup_{1\leq i,j\leq n}\beta(2\vartheta_{ij}+3)\leq \theta_4<\frac{1}{2}-\beta \sup_{1\leq i,j\leq n}\vartheta_{ij},\label{30}\\
&\sup_{1\leq i,j\leq n}\beta(\vartheta_{ij}+1)\leq \theta_5<\frac{1}{2}-\beta(\sup_{1\leq i,j\leq n}\vartheta_{ij}+1).\label{31}
\end{align}
By \eqref{26}, \eqref{28} and \eqref{29}, we derive that $\alpha$ should satisfy
\begin{align*}
\beta(3\sup_{1\leq i,j\leq n}\vartheta_{ij}+4)\leq \alpha<\frac{1}{2}-\beta(2\sup_{1\leq i,j\leq n}\vartheta_{ij}+2).
\end{align*}
Therefore, $\beta$ has to satisfy
\begin{align*}
0<\beta<\frac{1}{2(5\sup_{1\leq i,j\leq n}\vartheta_{ij}+6)}.
\end{align*}
We take $\kappa$ big enough to ensure that for any $\tilde \lambda>0$ it holds
\begin{align*}
\sup_{1\leq i,j\leq n}N^{\kappa[\beta(\vartheta_{ij}+2)+\alpha-\theta_2]}\leq N^{-\tilde\lambda},\quad
N^{\kappa(\alpha-\theta_3)}\leq N^{-\tilde \lambda}.
\end{align*}
Then we choose $m_1$ - $m_8$ big enough such that
\begin{align*}
&C\sup_{1\leq i,j\leq n}\big(N^{\kappa\beta(2\vartheta_{ij}+3)+2m_1(\theta_1-\frac{1}{2}+\beta \vartheta_{ij})+1}
+N^{\kappa\beta(\vartheta_{ij}+2)+2m_7[-\frac{1}{2}+\beta(\vartheta_{ij}+2)]+1}\\
&+N^{\kappa\beta(\vartheta_{ij}+2)+2m_8[-\frac{1}{2}+\beta(\vartheta_{ij}+2)]+1}
+N^{\kappa[\beta(2\vartheta_{ij}+2)+\alpha]+2m_2(\theta_2-\frac{1}{2}+\beta \vartheta_{ij})+1}\\
&+N^{\kappa[\beta(\vartheta_{ij}+1)+\alpha]+2m_3[\theta_3-\frac{1}{2}+\beta (\vartheta_{ij}+1)]+1}
+N^{\kappa\beta(3\vartheta_{ij}+3)+2m_4(\theta_4-\frac{1}{2}+\beta \vartheta_{ij})+1}\\
&+N^{\kappa\beta(2\vartheta_{ij}+2)+2m_5[\theta_5-\frac{1}{2}+\beta (\vartheta_{ij}+1)]+1}
+N^{\kappa\beta(\vartheta_{ij}+1)+2m_6[-\frac{1}{2}+\beta (\vartheta_{ij}+1)]+1}\big)\\
\leq& CN^{-\tilde\lambda}.
\end{align*}
As a consequence, we infer that
\begin{align*}
\mathbb{E}(S^\kappa_\alpha(t))\leq C\ln N\int_0^t \mathbb{E}(S^\kappa_\alpha(s))\,ds+CN^{-\tilde\lambda}.
\end{align*}
By means of the Gr\"onwall inequality, we deduce that
\begin{align*}
\sup_{0\leq t\leq T}\mathbb{P}\Big(\sum^n_{i=1}\max_{1\leq k\leq N}|X^{N,k}_{\eta,i}(t)-\bar X^k_{\eta,i}(t)|>N^{-\alpha}\Big)
\leq \sup_{0\leq t\leq T}\mathbb{E}(S^\kappa_\alpha(t))\leq CN^{-\tilde \lambda+C(T)}.
\end{align*}
Taking $\tilde \lambda=\lambda+C(T)$, it holds
\begin{align*}
\sup_{0\leq t\leq T}\mathbb{P}\Big(\sum^n_{i=1}\max_{1\leq k\leq N}|X^{N,k}_{\eta,i}(t)-\bar X^k_{\eta,i}(t)|>N^{-\alpha}\Big)
\leq \sup_{0\leq t\leq T}\mathbb{E}(S^\kappa_\alpha(t))\leq CN^{-\lambda}.
\end{align*}
For any $\tilde \alpha\leq \alpha$, we have
\begin{align*}
&\sup_{0\leq t\leq T}\mathbb{P}\Big(\sum^n_{i=1}\max_{1\leq k\leq N}|X^{N,k}_{\eta,i}(t)-\bar X^k_{\eta,i}(t)|>N^{-\tilde\alpha}\Big)\\
\leq&\sup_{0\leq t\leq T}\mathbb{P}\Big(\sum^n_{i=1}\max_{1\leq k\leq N}|X^{N,k}_{\eta,i}(t)-\bar X^k_{\eta,i}(t)|>N^{-\alpha}\Big)
\leq CN^{-\lambda}.
\end{align*}

\medskip
\indent
{\bf Acknowledgements:}
Yue Li would like to thank Chair in Applied Analysis of the University of Mannheim for
hosting her one-year scientific visit as an exchange doctoral student.
Y. Li are supported by NSFC (Grant No. 12071212).
Z. Zhang is supported by NSFC (Grant No. 12101305).
The authors thank the reviewer for her/his constructive comments and helpful suggestions.

\end{document}